\documentclass[number,citesort,seceqn,dvips]{arxbj}

\aid{0}
\volume{17}
\issue{3}
\pubyear{2011}
\firstpage{969}
\lastpage{986}
\doi{10.3150/10-BEJ306}

\makeatletter
\def\pr{\mathsf{P}}
\def\ep{\mathsf{E}}
\newtheorem{lemma}{Lemma}[section]

\newtheorem{theorem}{Theorem}[section]

\newremark{remark}{Remark}[section]
\makeatother

\begin{document}
\begin{frontmatter}

\title{On non-stationary threshold autoregressive models}
\runtitle{On non-stationary threshold autoregressive models}

\begin{aug}
\author{\fnms{Weidong} \snm{Liu}\thanksref{e1}\ead[label=e1,mark]{liuweidong99@gmail.com}},
\author{\fnms{Shiqing} \snm{Ling}\corref{}\thanksref{e2}\ead[label=e2,mark]{maling@ust.hk}}
\and
\author{\fnms{Qi-Man} \snm{Shao}\thanksref{e3}\ead[label=e3,mark]{maqmshao@ust.hk}}
\runauthor{W. Liu, S. Ling and Q.-M. Shao}
\address{Department of Mathematics,
Hong Kong University of Science and Technology,
Clear Water Bay, Kowloon,
Hong Kong, China.
\printead{e1}; \printead*{e2}; \printead*{e3}}
\end{aug}

\received{\smonth{5} \syear{2009}}
\revised{\smonth{6} \syear{2010}}

%
\begin{abstract}
In this paper we study the limiting
distributions of the least-squares estimators for the non-stationary
first-order threshold autoregressive ($\operatorname{TAR}(1)$) model.
It is proved
that the limiting behaviors of the $\operatorname{TAR}(1)$ process are very
different from those of the classical unit root model and the
explosive $\operatorname{AR}(1)$.
\end{abstract}

%

\end{frontmatter}
%
\section{Introduction}\label{sec1}

Since \cite{15}, threshold autoregressive (TAR) models have been
extensively investigated in the literature. The standard $\operatorname{TAR}(1)$
model can be written as follows:
%
%
\begin{equation}\label{mo1}
Y_{t}=\cases{
\gamma+\alpha Y_{t-1}+\varepsilon_{t}, & \quad $\mbox{if } Y_{t-1}>
r,$\vspace*{2pt}\cr
\delta+\beta Y_{t-1}+\varepsilon_{t}, &\quad $\mbox{if } Y_{t-1}\leq r,$}
\end{equation}
where $\{\varepsilon_{n}\}$ is a sequence of i.i.d. random
variables with zero mean and a finite variance $\sigma^{2}>0$.
Petrucceli and Woolford \cite{12} and Chan et al.
\cite{6} showed that, if $\varepsilon_{n}$ has a strictly
positive density, then the necessary and sufficient condition for
the strictly stationary and geometrically ergodic solution to model
(\ref{mo1}) when $\gamma=\delta=0$ is
%
%
\begin{equation}\label{cd1}
\alpha<1,\qquad  \beta<1 \quad  \mbox{and} \quad \alpha\beta<1;
\end{equation}
see also \cite{14,5}. The
properties of the least-squares estimator (LSE) of model (\ref{mo1})
were established when $\{Y_{t}\}$ is stationary by Chan \cite{4} and
later by Chan and Tsay \cite{5} for the continuous case (i.e.,
$\gamma+r\alpha=\delta+r\beta$). When $(\alpha, \beta)$ does not lie
in the stationary~re\-gion~(\ref{cd1}), the estimation theory of the LSE of
model (\ref{mo1}) is challenging.

Pham, Chan and Tong \cite{13} were the first to
consider the non-stationary case of~mo\-del~(\ref{mo1}). They focus
on the following case:
%
%
\begin{equation}\label{ca1}
\gamma=\delta   \quad \mbox{and} \quad      r=0
\end{equation}
and assume that $\delta$ is a known parameter. For the
LSE of $(\alpha,\beta)$:
%
%
\begin{eqnarray}\label{de1}
\hat{\alpha}_{n}&=&\frac{\sum_{t=1}^{n-1}I(Y_{t}>
r)Y_{t}(Y_{t+1}-\gamma)}{\sum_{t=1}^{n-1}I(Y_{t}> r)Y^{2}_{t}},
\\
 \hat{\beta}_{n}&=& \frac{\sum_{t=1}^{n-1}I(Y_{t}\leq
r)Y_{t}(Y_{t+1}-\gamma)}{\sum_{t=1}^{n-1}I(Y_{t}\leq
r)Y^{2}_{t}},\nonumber
\end{eqnarray}
they show that
\[
(\hat{\alpha}_{n},\hat{\beta}_{n})\rightarrow
(\alpha,\beta)\qquad    \mbox{a.s.}
\]
if and only if one of the following conditions holds:
%
%
\begin{eqnarray}\label{eq1.5}
\alpha&\leq&1, \qquad  \beta\leq1\quad  \mbox{and}\quad  \gamma=0,\nonumber\\
\alpha&<&1, \qquad \beta\leq1\quad   \mbox{and}\quad  \gamma>0,\\
\alpha&\leq&1, \qquad  \beta<1\quad   \mbox{and} \quad
\gamma<0.\nonumber
\end{eqnarray}
They also showed that, when $\alpha\beta=1$, the estimator of
$\alpha$ is strongly consistent. However, the rate of convergence
and the limiting distribution of LSE are two open questions when
$(\alpha, \beta)$ lies in the non-ergodic region.

Following \cite{13}, in this paper we study the limiting
distribution of $(\hat{\alpha}_{n}, \hat{\beta}_{n})$ for the
following cases:
\begin{eqnarray*}
&&\mbox{Case I:}\quad  \gamma=\delta=0,\qquad  \alpha= 1\quad \mbox{and} \quad \beta< 1,\\
&&\mbox{Case II:}\quad  \gamma=\delta=0,\qquad \alpha>1 \quad\mbox{and} \quad \beta\leq1.
\end{eqnarray*}
For each case, we partially derive the limiting distribution of
$(\hat{\alpha}_{n}, \hat{\beta}_{n})$ under some suitable
conditions. Case I is related to the unit root problem, which is
particularly interesting in economics and finance. One usually tests
whether or not a market is efficient via testing a~unit root in AR
model. Unit root tests have been extensively studied in the
literature; see, for example, \cite{9,3,10}. When $Y_{t}$ denotes a market index,
case I can describe the phenomena that the market moves from
efficiency to inefficiency when the index crosses the threshold $r$ and
$|\beta|<1$.
Our result may provide a way to test this phenomena. The results for
case II can help
us to understand the limiting behaviors of the LSE in this complicated
and dynamic system.
Our proof is based on the limiting behavior of~$Y_{t}$ as
$t\to\infty$. The method of the proof is non-standard and may
provide some insights for future research in this area.

The paper is organized as follows. The main results are stated in
Section \ref{sec2}. The proofs of the main results are given in Sections
\ref{sec3} and \ref{sec4}. This paper also includes consistency of the LSE when
$\alpha\beta=1$ in Section \ref{sec5}, which is of independent interest.
Throughout, we let $C$ and~$C_{(\cdot)}$ denote positive constants
that may be different in every place, $\mathcal{
F}_{t}=\sigma\{Y_{0},\varepsilon_{1},\ldots, \varepsilon_{t}\}$, and
we assume the initial value $Y_{0}$ in model (\ref{mo1}) is a random
variable independent with $n$ and $\{\varepsilon_{t}; t\geq1\}$.

\section{Main results}\label{sec2}

We consider two different cases.

\subsection{Case I}

Assume $r\le0$ or $r> 0$ with
$\alpha=1$ and $\beta=-1$. The results are stated as follows.

\begin{theorem}\label{th1} Assume $\gamma=\delta=0$ and $\ep
Y^{2}_{0}<\infty$. If either
\begin{longlist}[(ii)]
\item[(i)] $\alpha=1$, $\beta< 1$ and $r\leq0$; or
\item[(ii)] $\alpha=1$, $\beta=-1$ and $r\in R$;
is satisfied, then we have
%
%
\begin{equation}\label{p2}
n(\hat{\alpha}_{n}-1)\Rightarrow
\frac{B^{2}(1)-1}{2\int_{0}^{1}B^{2}(t)\,\mathrm{d}t},
\end{equation}
where $B(t)$ is a standard Brownian motion.
\end{longlist}
\end{theorem}

\begin{remark} Unlike the stationary case in \cite{4}, the limiting distribution of $\hat{\alpha}_{n}$ is
independent of~$r$ and $\beta$. (\ref{p2}) could be used to test
whether $(\alpha,\beta)$ lies on the boundary $\{\alpha=1,\beta<1\}$
if we know~$r$ is zero or negative. We note that this test is the
same as the Dickey--Fuller test. The limiting distribution of
$\hat{\beta}_{n}$ is still unclear. But when $\alpha<1$, $\beta=1$
and $r\geq0$, from the proof of Theorem \ref{th1} and Remark \ref{rem3.1},
we have
%
%
\begin{equation}\label{p2.b}
n(\hat{\beta}_{n}-1)\Rightarrow
\frac{B^{2}(1)-1}{2\int_{0}^{1}B^{2}(t)\,\mathrm{d}t}.
\end{equation}
We should mention that Caner and Hansen \cite{7} developed an asymptotic
theory for a~TAR model with a unit root, but
their model is not the same as model (\ref{mo1}) since their threshold
variable is $Y_{t-1}-Y_{t-2}$.
\end{remark}

\begin{remark} When $\alpha=1$, $\beta<1$ and
$\gamma=\delta<0$, Chan \textit{et al.} \cite{6} show that $\{Y_{t}\}$ is
ergodic, and hence is strictly stationary by assuming that $Y_{0}$
has its distribution $\pi(\cdot)$ that is the invariant probability
distribution of $\{Y_{t}\}$. For the case $\alpha=1$, $\beta<1$,
$r\leq0$ and $\gamma=\delta>0$, we have
\[
Y_{n}\geq
\gamma+\varepsilon_{n}+Y_{n-1}\geq
n\gamma+\sum_{k=1}^{n}\varepsilon_{k}+Y_{0}.
\]
Hence
$Y_{n}\rightarrow\infty$ a.s. It follows that $\max_{1\leq k\leq
n}|Y_{k}-k\gamma-\sum_{i=1}^{k}\varepsilon_{i}|=\mathrm{O}(1)$ a.s. By some
standard arguments using the martingale central limit theorem (CLT), it
is not
hard to see $n^{3/2}(\hat{\alpha}_{n}-1)\Rightarrow
N(0,3\sigma^{2}/\gamma^{2})$. In this case, $\hat{\beta}_{n}$ is not
a strongly consistent estimator.
\end{remark}

\subsection{Case II}
By (\ref{eq1.5}), $(\hat{\alpha}_{n},\hat{\beta}_{n})$ is not a consistent
estimator of $(\alpha,\beta)$ in this case. However, the following
theorem shows that $\hat{\alpha}_{n}$ is a consistent estimator of
$\alpha$.

\begin{theorem}\label{th3} Assume that $\gamma=\delta=0$ and one of the
following conditions holds:
\begin{longlist}[(H1)]
\item[(H1)] $\alpha>1$, $\beta\leq1$, $r=0$ and $\ep
Y^{2}_{0}<\infty$;

\item[(H2)] $\alpha>1$, $\beta\le1$, $r\not=0$, $\ep
Y_{0}^{2}<\infty$ and $\pr(\varepsilon_{1}\leq x)<1$ for any $x\in
R$.
\end{longlist}
Then we have
\[
(\alpha^{2}-1)^{-1}\alpha^{n}(\hat{\alpha}_{n}-\alpha)\Rightarrow
\eta^{*}/\xi^{*},
\]
where $\eta^{*}$ and $\xi^{*}$ are independent
random variables, $\eta^{*}\stackrel{d}{=}
\sum_{t=1}^{\infty}\alpha^{-t}\varepsilon_{t}$,
$\xi^{*}\stackrel{d}{=} \xi$ and
%
%
\begin{eqnarray}\label{newadd-11}
\xi&=&\sum_{k=1}^{\infty}\alpha^{-k+1}\biggl(\frac{ \beta
}{\alpha}\biggr)^{m_{k}}\varepsilon_{k}+\biggl(\frac{\beta}{\alpha
}
\biggr)^{m_{0}}Y_{0}>0\qquad
\mbox{a.s. for $\beta\leq1$ and $\beta\neq0$;}\\
\label{eq2.4}\xi&=&\sum_{k=1}^{\infty}\frac{\prod_{t=k}^{\infty}I\{Y_{t}>
r\}}{\alpha^{k}}\varepsilon_{k} +\prod_{t=0}^{\infty}I\{Y_{t}>
r\}Y_{0}>0\qquad   \mbox{a.s. for $\beta=0$},
\end{eqnarray}
where $m_{k}=\sum_{t=k}^{\infty}I\{Y_{t}\le r\}$ is almost surely
finite.
\end{theorem}

\begin{remark}
In the explosive $\operatorname{AR}(1)$ model, $Y_{t}=\alpha Y_{t-1}+\varepsilon_{t}$,
it is well known that the LSE of~$\alpha$ asymptotically follows
a Cauchy distribution if $\varepsilon_{t}$ is normal. By Theorem \ref{th3},
this conclusion does not hold any more for model (\ref{mo1}).
\end{remark}

\section{\texorpdfstring{Proof of Theorem \protect\ref{th1}}{Proof of Theorem 2.1}}\label{sec3}

Before proving Theorem \ref{th1}, we first establish the limiting
distribution for $\{Y_{t}\}$ when $t\rightarrow\infty$ as follows.

\begin{theorem}\label{th4.1} If either \textup{(i)} or \textup{(ii)} in Theorem \ref
{th1} holds,
then
%
%
\begin{equation}\label{p1}\frac{Y_{[nt]}}{\sqrt{n}}\Rightarrow
\sigma|B(t)|\qquad \mbox{on $D[0,1]$},
\end{equation}
as $n\to\infty$, where $D[0,1]$ is the Skorokhod space.
\end{theorem}

\begin{remark}\label{rem3.1} It is interesting to see that the
limiting distribution in (\ref{p1}) does not depend on $\beta$ and~$r$. This means that the effect of $\beta$ and $r$ on $Y_{t}$ is
ignorable when $t$ is long enough. The pattern of $Y_{t}$ is quite
different from the unit root process in the $\operatorname{AR}(1)$ model in which
$X_{[nt]}/{\sqrt{n}}\Rightarrow\sigma B(t)$ on $D[0,1]$, where
$X_{t}=X_{t-1}+\varepsilon_{t}$. If $\beta=1$, $\alpha<1$ and $r\geq
0$, then replacing $Y_{t}$ by $-Y_{t}$, we can get
$Y_{[nt]}/\sqrt{n}\Rightarrow-\sigma|B(t)|$ on $D[0,1]$.
\end{remark}

\begin{pf*}{Proof of Theorem \protect\ref{th4.1} under (ii)}  We
first consider the case when $\alpha=1$ and $\beta=-1$. Denote~$Y_{n}$ by $Y^{\star}_{n}$ in this case. If $r\geq0$, we have
$Y^{\star}_{n}=\varepsilon_{n}+|Y^{\star}_{n-1}|-2Y^{\star}_{n-1}I\{
0\leq
Y^{\star}_{n-1}\leq r\}$, and if $r<0$, we have
$Y^{\star}_{n}=\varepsilon_{n}+|Y^{\star}_{n-1}|+2Y^{\star}_{n-1}I\{r<
Y^{\star}_{n-1}\leq0\}$. Hence,
\[
\max_{1\leq k\leq
n}\bigl|Y^{\star}_{k}-|Y^{\star}_{k-1}|\bigr|\leq\max_{1\leq k\leq
n}|\varepsilon_{k}|+2|r|=\mathrm{o}_{\pr}\bigl(\sqrt{n}\bigr).
\]
So it is enough to
show that $|Y^{\star}_{[nt]}|/\sqrt{n}\Rightarrow\sigma|B(t)|$ on
$D[0,1]$. Note that
\[
Y^{\star}_{n}=\sum_{k=1}^{n}\prod_{j=k}^{n-1}I_{j}\varepsilon
_{k}+\prod
_{j=0}^{n-1}I_{j}Y^{\star}_{0},
\]
where $I_{k}=I\{Y^{\star}_{k}>r\}-I\{Y^{\star}_{k}\leq r\}$. It
follows that
\[
A^{-1}_{n}Y^{\star}_{n}=\sum_{k=1}^{n}A^{-1}_{k}\varepsilon
_{k}+Y^{\star}_{0},
\]
where $A_{n}=\prod_{k=0}^{n-1}I_{k}$. Since $\ep
[A^{-2}_{k}\varepsilon^{2}_{k}|\mathcal{F}_{k-1}]=1$, we have by the
martingale CLT (cf. \cite{2}) that
%
%
\begin{eqnarray}
\frac{1}{\sqrt{n}}A^{-1}_{[nt]}Y^{\star}_{[nt]}\Rightarrow\sigma
B(t).
\end{eqnarray}
Now, (\ref{p1}) follows from $|A_{k}|=1$ and the continuous
mapping theorem.
\end{pf*}

\begin{pf*}{Proof of Theorem \protect\ref{th4.1} under (i)}
Recall
the definition of $\{Y^{\star}_{n}\}$ with the initial value
$Y^{\star}_{0}=Y_{0}$. For any $p>0$, observe that
%
\begin{eqnarray}\label{aeq17}
|Y_{n}-Y^{\star}_{n}|^{p}&=&|Y_{n-1}-Y^{\star}_{n-1}|^{p}I\{
Y_{n-1}>r,Y^{\star}_{n-1}>r\}\nonumber\\
& &{}+|Y_{n-1}+Y^{\star}_{n-1}|^{p}I\{Y_{n-1}>r,Y^{\star}_{n-1}\leq
r\}
\nonumber
\\[-8pt]
\\[-8pt]
\nonumber
& &{}+|\beta Y_{n-1}-Y^{\star}_{n-1}|^{p}I\{Y_{n-1}\leq
r,Y^{\star}_{n-1}> r\}\\
& &{}+|\beta
Y_{n-1}+Y^{\star}_{n-1}|^{p}I\{Y_{n-1}\leq r,Y^{\star}_{n-1}\leq
r\}.\nonumber
\end{eqnarray}
Since $r\leq0$, it follows that
%
\begin{eqnarray}\label{aeq18}
&&|Y_{n-1}+Y^{\star}_{n-1}|^{p}I\{Y_{n-1}>r,Y^{\star}_{n-1}\leq
r\}\nonumber\\
&&\quad   \leq
|Y_{n-1}-Y^{\star}_{n-1}|^{p}I\{Y_{n-1}>0,Y^{\star}_{n-1}\leq
r\}\\
&&\qquad{}   +C_{p}|Y^{\star}_{n-1}|^{p}I\{Y^{\star}_{n-1}\leq
r\}+C_{p,r}I\{Y_{n-1}\leq0\}.\nonumber
\end{eqnarray}
Furthermore, since $\beta\leq1$ and $r\leq0$, we have
%
\begin{eqnarray}\label{aeq19}
&&|\beta Y_{n-1}-Y^{\star}_{n-1}|^{p}I\{Y_{n-1}\leq
r,Y^{\star}_{n-1}> r,Y^{\star}_{n-1}\geq\beta Y_{n-1}\}\nonumber\\
&&\quad  \leq
|Y_{n-1}-Y^{\star}_{n-1}|^{p}I\{Y_{n-1}\leq r,Y^{\star}_{n-1}>
r,Y^{\star}_{n-1}\geq\beta Y_{n-1}\},
\nonumber
\\[-8pt]
\\[-8pt]
\nonumber
&&|\beta
Y_{n-1}-Y^{\star}_{n-1}|^{p}I\{Y_{n-1}\leq r,Y^{\star}_{n-1}>
r,Y^{\star}_{n-1}< \beta Y_{n-1,\beta}\}\\
&&\quad
\leq
2^{p}|\beta Y_{n-1}|^{p}I\{Y_{n-1}\leq
r\}+C_{p,\beta,r}I\{Y_{n-1}\leq r\}.\nonumber
\end{eqnarray}
It follows from (\ref{aeq17})--(\ref{aeq19}) that
%
%
\begin{equation}\label{eq10}
|Y_{n}-Y^{\star}_{n}|^{p}\leq|Y_{n-1}-Y^{\star}_{n-1}|^{p}+q_{n}\le
\sum_{k=1}^{n}q_{k},
\end{equation}
where
%
%
\begin{equation}\label{eq11}
q_{n}=C|Y_{n-1}|^{p}I\{Y_{n-1}\leq
r\}+C|Y^{\star}_{n-1}|^{p}I\{Y^{\star}_{n-1}\leq r\}+CI\{Y_{n-1}\leq
0\}.
\end{equation}

We first have $\pr(Y_{n}\leq r)\rightarrow0$ by Lemma
\ref{le6.2} below and $\pr(Y^{\star}_{n-1}\leq r)\to0$ by
(\ref{p1}) under $\alpha=1$ and $\beta=-1$ as $n\rightarrow\infty$.
Furthermore, applying Lemma \ref{le3.1} below with $p=2$ for $Y_{t}$
and $Y_{t}^{\star}$, we have
\begin{eqnarray*}
\ep q_{n}&\leq&
C\sup_{k}\ep(\varepsilon^{2}_{k}+Y^{2}_{0})I\{Y_{n-1}\leq
r\}+C\sup_{k}\ep(\varepsilon^{2}_{k}+Y^{\star
2}_{0})I\{Y^{\star}_{n-1}\leq
r\}\\
& &{} +C\pr(Y_{n-1}\leq r)+C\pr(Y^{\star
}_{n-1}\leq
r)\\
&\rightarrow& 0.
\end{eqnarray*}
Thus, by (\ref{aeq17}) with $p=2$ and the previous inequality, we have
\[
\ep\max_{1\leq k\leq
n}{|}Y_{k}-Y^{\star}_{k}{|}^{2}\big/n\le
\frac{1}{n}\sum_{k=1}^{n}\ep q_{k}\rightarrow0
\]
as $n\rightarrow\infty$. By (i) of Theorem \ref{th4.1} and the previous
inequality, (ii) of Theorem \ref{th4.1} holds. This completes the proof.
\end{pf*}

We are now ready to prove Theorem \ref{th1}.

\begin{pf*}{Proof of Theorem \protect\ref{th1}}
Without loss of
generality, we assume $\sigma=1$. Note that
\[
n(\hat{\alpha}_{n}-1)=\frac{n\sum_{t=1}^{n-1}I\{Y_{t}>r\}
Y_{t}\varepsilon_{t+1}}{\sum_{t=1}^{n-1}I\{Y_{t}>
r\}Y^2_{t}}
\]
and
\begin{eqnarray*}
\frac{1}{n}\sum_{t=1}^{n-1}I\{Y_{t}>r\}Y_{t}\varepsilon_{t+1}
&=&\frac{1}{2n}\sum_{t=1}^{n-1}I\{Y_{t}>r\}
[(Y^{2}_{t+1}-Y^{2}_{t})-\varepsilon^{2}_{t+1}]\\
&=&\frac{1}{2n}\sum_{t=1}^{n-1}(Y^{2}_{t+1}-Y^{2}_{t})-\frac
{1}{2n}\sum
_{t=1}^{n-1}I\{Y_{t}\leq
r\}(Y^{2}_{t+1}-Y^{2}_{t})\\
&&{} -\frac{1}{2n}\sum_{t=1}^{n-1}I\{Y_{t}>r\}
\varepsilon^{2}_{t+1}.
\end{eqnarray*}
Since $\pr(y_{n}\leq r)\rightarrow\pr(|B(1)|\leq0)=0$ as
$n\to\infty$, we have
\[
\frac{1}{n}\sum_{t=1}^{n-1}\ep I\{Y_{t}\leq r\}
\varepsilon^{2}_{t+1}=\frac{1}{n}\sum_{t=1}^{n-1}\pr(Y_{t}\leq
r)\rightarrow0,
\]
as $n\to\infty$. Thus,
%
%
\begin{equation}\label{eqnew-1}
\frac{1}{n}\sum_{t=1}^{n-1}I\{Y_{t}>r\}
\varepsilon^{2}_{t+1}\rightarrow1
\end{equation}
in probability. Furthermore, we
have $n^{-1}\sum_{t=1}^{n-1}\ep Y^{2}_{t}I\{Y_{t}\leq r\}\rightarrow
0$ by Lemmas~\ref{le3.1} and~\ref{le6.2} below, and hence
%
%
\begin{equation}\label{eq2}
\frac{1}{n}\sum_{t=1}^{n-1}I\{Y_{t}\leq
r\}(Y^{2}_{t+1}-Y^{2}_{t})\rightarrow0
\end{equation}
in probability. Thus, by (\ref{eqnew-1}) and (\ref{eq2}), we have
\[
\frac{1}{n}\sum_{t=1}^{n-1}I\{Y_{t}>r\}Y_{t}\varepsilon_{t+1}
=\frac{1}{2n}(Y^{2}_{n}-Y^{2}_{1})-\frac{1}{2}+\mathrm{o}_{\pr}(1)\Rightarrow
\frac{1}{2}B^{2}(1)-\frac{1}{2}.
\]
Note that
\[
\frac{\sum_{t=1}^{n-1}I\{Y_{t}>
r\}Y^2_{t}}{n^2}=\int_{0}^{1}I\bigl\{Y_{[nt]}>r\bigr\}\frac
{Y^{2}_{[nt]}}{n}\,\mathrm{d}t\Rightarrow
\int_{0}^{1}B^{2}(t)\,\mathrm{d}t.
\]
Theorem \ref{th1} follows from the continuous mapping theorem.
\end{pf*}

We now prove Lemmas \ref{le3.1} and \ref{le6.2}, which were used in the
proof of Theorem \ref{th1}.

\begin{lemma}\label{le3.1} Suppose that $\ep|Y_{0}|^{p}<\infty$ and
$\ep|\varepsilon_{0}|^{p}<\infty$ for some $p>0$.
Under the conditions $\gamma=\delta=0$, $\alpha=1$ and $\beta<1$,
for any event $A$, it holds that
\[
\ep|Y_{n}|^{p}I\{Y_{n}\leq r,A\}\leq
C\Bigl(\sup_{k}\ep|\varepsilon_{k}|^{p}I\{A\}+\ep
|Y_{0}|^{p}I\{A\}+\pr(A)\Bigr).
\]
\end{lemma}

\begin{pf}
Set $X_{n}=Y_{n}-r$ for $n\geq0$. We can see that
$X_{n}=e_{n}+X^{+}_{n-1}-\beta X_{n-1}^{-}$, where
$e_{n}=\varepsilon_{n}+(\beta-1)rI\{X_{n-1}\leq0\}$. Suppose
$\beta\leq0$. The lemma follows from
%
%
\begin{eqnarray}\label{eq9}
&&\ep|X_{n}|^{p}I\{X_{n}\leq0,A\}\nonumber\\
&&\quad =\ep|X_{n}|^{p}I\{e_{n}\leq
-(X^{+}_{n-1},-\beta X_{n-1}^{-}),A\}
\nonumber
\\[-8pt]
\\[-8pt]
\nonumber
&&\quad\leq
C_{p}\ep|e_{n}|^{p}I\{A\}+C_{p}\ep|X_{n-1}^{+}-\beta
X_{n-1}^{-}|^{p}I\{|e_{n}|\geq X_{n-1}^{+}-\beta X_{n-1}^{-},A\}\\
&&\quad\leq
2C_{p}\sup_{k}\ep|\varepsilon_{k}|^{p}I\{A\}+C_{p,\beta,r}\pr(A).\nonumber
\end{eqnarray}
Now we prove the lemma when $0<\beta<1$. Set the events
$A_{k}=\{X_{k}\leq0\}$ for $1\leq k\leq n$. Note that
%
\begin{eqnarray}\label{p5}\nonumber
\ep|X_{n}|^{p}I\{X_{n}\leq0,A\}&=&\sum_{k=0}^{n-1}\ep|X_{n}|^{p}I\{A_{n}\cdots A_{n-k}A^{c}_{n-k-1},A\}\\
&&+\ep|X_{n}|^{p}I\{A_{n}\cdots A_{0},A\}.
\end{eqnarray}
We need to estimate $\ep|X_{n}|^{p}I\{A_{n}\cdots
A_{n-k}A^{c}_{n-k-1}\}$. In fact, on $A_{n}\cdots
A_{n-k}A^{c}_{n-k-1}$, we have
\[
X_{n}=\sum_{j=0}^{k-1}\beta^{j}e_{n-j}+\beta^{k}X_{n-k}.
\]
Set $\xi=\sum_{j=0}^{n}\beta^{j}|e_{n-j}|$. It follows that
%
\begin{eqnarray}\label{eq3}
&&\ep|X_{n}|^{p}I\{A_{n}\cdots A_{n-k}A^{c}_{n-k-1},A\}\nonumber\\
&&\quad \leq
C_{\beta,\delta}\ep|\xi|^{p}I\{A_{n}\cdots
A_{n-k}A^{c}_{n-k-1},A\}\nonumber\\
&&\qquad{}+C_{p}\beta^{kp}\ep
|X_{n-k}|^{p}I\{A_{n}\cdots A_{n-k}A^{c}_{n-k-1},A\}\nonumber\\
&&\quad \leq
C_{p,\beta}\ep|\xi|^{p}I\{A_{n}\cdots
A_{n-k}A^{c}_{n-k-1},A\}
\nonumber
\\[-8pt]
\\[-8pt]
\nonumber
&&\qquad{}+C_{p}\beta^{kp}\Bigl(\sup_{k}\ep
|\varepsilon_{k}|^{p}I\{A\}+C_{p,\beta,r}\pr(A)\Bigr)\\
&&\qquad{}
+C_{p}\beta^{kp}\ep|X_{n-k-1}|^{p}I\{|e_{n-k}|\geq
X_{n-k-1},A^{c}_{n-k-1},A\}\nonumber\\
&&\quad \leq
C_{p,\beta}\ep|\xi|^{p}I\{A_{n}\cdots
A_{n-k}A^{c}_{n-k-1},A\}\nonumber\\
&&\qquad{}+2C_{p}\beta^{kp}\Bigl(\sup_{k}\ep
|\varepsilon_{k}|^{p}I\{A\}+C_{p,\beta,r}\pr(A)\Bigr).\nonumber
\end{eqnarray}
Clearly, on $A_{n}\cdots A_{0}$, we have
$X_{n}=\sum_{j=0}^{n-1}\beta^{j}e_{n-j}+\beta^{n}X_{0}$ and hence by
(\ref{p5}) and~(\ref{eq3}), $\ep|X_{n}|^{p}I\{X_{n}\leq0,A\}\leq
C(\sup_{k}\ep|\varepsilon_{k}|^{p}I\{A\}+\ep
|Y_{0}|^{p}I\{A\}+\pr(A))$. The lemma is now proved.
\end{pf}

\begin{lemma}\label{le6.2} Suppose that $\ep Y_{0}^{2}<\infty$, $\ep
\varepsilon_{0}=0$ and $\ep\varepsilon_{0}^{2}<\infty$.
Under the conditions $\gamma=\delta=0$, $\alpha=1$, $\beta<1$ and
$r\leq0$, we have $Y_{n}/\sqrt{n}\Rightarrow\sigma|B(1)|$
as $n\rightarrow\infty$.
\end{lemma}
\begin{pf}
For $K>0$, set
\begin{eqnarray*}
\widetilde{\varepsilon}_{k}&=&\varepsilon_{k}I\{|\varepsilon
_{k}|\leq
K\}-\ep\varepsilon_{k}I\{|\varepsilon_{k}|\leq K\},\qquad
\hat{\varepsilon}_{k}=\varepsilon_{k}-\widetilde{\varepsilon
}_{k},\qquad
k\geq1.\\
 \widetilde{Y}_{0}&=&Y_{0}I\{|Y_{0}|\leq K\},\qquad
\widehat{Y}_{0}=Y_{0}-\widetilde{Y}_{0}.
\end{eqnarray*}
We now construct two $\operatorname{TAR}(1)$ processes $\{\widetilde{Y}_{t}\}$ and
$\{\widetilde{Y}_{t}^{\star}\}$ as follows:
%
%
\begin{eqnarray}\label{6.4}
\widetilde{Y}_{n}&=&\widetilde{\varepsilon}_{n}+\widetilde{Y}_{n-1}I\{
\widetilde{Y}_{n-1}>r\}+
\beta\widetilde{Y}_{n-1}I\{\widetilde{Y}_{n-1}\leq r\},  \qquad n\geq
1;
\\
\label{6.44}
\widetilde{Y}^{\star}_{n}&=&\widetilde{\varepsilon}_{n}+\widetilde
{Y}^{\star}_{n-1}I\{\widetilde{Y}^{\star}_{n-1}>r\}
-\widetilde{Y}^{\star}_{n-1}I\{\widetilde{Y}^{\star}_{n-1}\leq r\},
 \qquad n\geq1.
\end{eqnarray}
By Theorem \ref{th4.1}, when $\alpha=1$ and $\beta=-1$, we can see
that
%
%
\begin{equation}\label{6.1}
\widetilde{Y}^{\star}_{[nt]}/\sqrt{n}\Rightarrow
\sigma_{K}|B(t)|  \qquad \mbox{on }   D[0,1],
\end{equation}
with
$\sigma^{2}_{K}=\operatorname{Var}(\widetilde{\varepsilon}_{1})$.
Let $q{}^{\prime}_{k}$,
$1\leq k\leq n$, be defined as $q_{k}$ in (\ref{eq11}) by replacing
$\{Y_{n}\}$ and $\{Y^{\star}_{n}\}$ with $\{\widetilde{Y}_{n}\}$ and
$\{\widetilde{Y}^{\star}_{n}\}$, respectively. Taking $p>2$, by
virtue of (\ref{eq10}), we have
\[
\ep\max_{1\leq k\leq
n}\biggl|\frac{\widetilde{Y}_{k}-\widetilde{Y}^{\star}_{k}}{\sqrt
{n}}\biggr|^{p}\leq
\frac{\sum_{k=1}^{n}\ep q{}^{\prime}_{k}}{n^{p/2}}.
\]
Furthermore, using Lemma \ref{le3.1} with $Y_{t}$ replaced by
$\{\widetilde{Y}_{n}\}$ and $\{\widetilde{Y}^{\star}_{n}\}$,
respectively, we know that $q{}^{\prime}_{k}$ is uniformly bounded for
all $k\ge1$. Thus, we have
\[
\ep\max_{1\leq k\leq
n}\biggl|\frac{\widetilde{Y}_{k}-\widetilde{Y}^{\star}_{k}}{\sqrt
{n}}\biggr|^{p}\leq
Cn^{-p/2+1}.
\]
By (\ref{6.1})
and the previous inequality, we have
%
%
\begin{equation}
\widetilde{Y}_{[nt]}/\sqrt{n}\Rightarrow\sigma_{K}|B(t)|
\qquad \mbox{on }  D[0,1].
\end{equation}
Since $\sigma_{K}\to\sigma$ as $K\to\infty$, it suffices to show
that for any $\delta>0$,
%
%
\begin{equation}\label{p3}
\lim_{K\rightarrow\infty}\limsup_{n\rightarrow\infty}\pr
\bigl(|Y_{n}-\widetilde{Y}_{n}|\geq
\delta\sqrt{n}\bigr)=0.
\end{equation}

By model (\ref{mo1}) and model (\ref{6.4}), we have
\begin{eqnarray*}
\ep(Y_{n}-\widetilde{Y}_{n})^{2}&=&\ep
\hat{\varepsilon}^{2}_{n}+\ep
(Y_{n-1}-\widetilde{Y}_{n-1})^{2}I\{Y_{n-1}>r,\widetilde
{Y}_{n-1}>r\}\\
& &{}+\ep(Y_{n-1}-\beta
\widetilde{Y}_{n-1})^{2}I\{Y_{n-1}>r,\widetilde{Y}_{n-1}\leq
r\}\\
& &{}+\ep(\beta
Y_{n-1}-\widetilde{Y}_{n-1})^{2}I\{Y_{n-1}\leq
r,\widetilde{Y}_{n-1}>r\}\\
& &{}+\ep(\beta
Y_{n-1}-\beta\widetilde{Y}_{n-1})^{2}I\{Y_{n-1}\leq
r,\widetilde{Y}_{n-1}\leq r\}.
\end{eqnarray*}
It can be verified that
\begin{eqnarray*}
&&\ep(Y_{n-1}-\beta
\widetilde{Y}_{n-1})^{2}I\{Y_{n-1}>r,\widetilde{Y}_{n-1}\leq
r\}\\
&&\quad  \leq\ep(Y_{n-1}-
\widetilde{Y}_{n-1})^{2}I\{Y_{n-1}>r,\widetilde{Y}_{n-1}\leq
r,Y_{n-1}>\beta\widetilde{Y}_{n-1}\}\\
&&\qquad{}  +C\ep
\widetilde{Y}^{2}_{n-1}I\{\widetilde{Y}_{n-1}\leq
r\}+C\pr(\widetilde{Y}_{n-1}\leq r).
\end{eqnarray*}
Let $M$ be any positive number. Then,
\begin{eqnarray*}
&&\ep(\beta
Y_{n-1}-\widetilde{Y}_{n-1})^{2}I\{Y_{n-1}\leq
r,\widetilde{Y}_{n-1}>r\}\nonumber\\
&&\quad \leq\ep(
Y_{n-1}-\widetilde{Y}_{n-1})^{2}I\{Y_{n-1}\leq
r,\widetilde{Y}_{n-1}>r,\widetilde{Y}_{n-1}>\beta Y_{n-1}\}\nonumber
\\
&&\qquad{}   +C\ep Y^{2}_{n-1}I\{Y_{n-1}<-M,\beta
Y_{n-1}\geq\widetilde{Y}_{n-1}>
r\}+C_{\beta,M,r}\pr(\widetilde{Y}_{n-1}\leq
r+|\beta|M).
\end{eqnarray*}
Combining the above inequalities, we can see that
\[
\ep(Y_{n}-\widetilde{Y}_{n})^{2}\leq\ep
\hat{\varepsilon}^{2}_{n}+\ep(Y_{n-1}-\widetilde
{Y}_{n-1})^{2}
+\widetilde{q}_{n},
\]
where
\begin{eqnarray*}
\widetilde{q}_{n}&=&C\ep Y^{2}_{n-1}I\{Y_{n-1}\leq
r,\widetilde{Y}_{n-1}\leq r\}+C\ep
\widetilde{Y}^{2}_{n-1}I\{\widetilde{Y}_{n-1}\leq r\}\\
& &{}+C\ep
Y^{2}_{n-1}I\{Y_{n-1}<-M,\beta Y_{n-1}\geq\widetilde{Y}_{n-1}>
r\}+C_{\beta,M,r}\pr(\widetilde{Y}_{n-1}\leq
r+|\beta|M).
\end{eqnarray*}
By induction we have
%
%
\begin{equation}\label{xd1}
\ep(Y_{n}-\widetilde{Y}_{n})^{2}\leq n\ep
\hat{\varepsilon}^{2}_{0}+\ep
\widehat{Y}^{\hspace*{2pt}2}_{0}+\sum_{k=1}^{n}\widetilde{q}_{k}.
\end{equation}
Since
$\widetilde{Y}_{n}/\sqrt{n}\Rightarrow\sigma_{K}|B(1)|$, we have
$\pr(\widetilde{Y}_{n-1}\leq r+|\beta|M)\rightarrow0$
as $n\rightarrow\infty$. Note that
\[
I\{y_{n-1}<-M,\beta Y_{n-1}\geq\widetilde{Y}_{n-1}>r\}\leq\cases{
I\{\varepsilon_{n-1}<-M+|r|\},&\quad $\mbox{if
$\beta\leq0$},$\vspace*{2pt}\cr
I\{\widetilde{Y}_{n-1}\leq-\beta M\},&\quad $\mbox{if $0<\beta<1$}.$}
\]
By Lemma \ref{le3.1},
\begin{eqnarray*}
&& \ep Y^{2}_{n-1}I\{Y_{n-1}<-M,\beta
Y_{n-1}\geq\widetilde{Y}_{n-1}> r\}\\
&&\quad  \leq
C\sup_{k}\ep(\varepsilon_{k}^{2}+Y^{2}_{0}+1)I\{\varepsilon
_{n-1}<-M+|r|\}\cr
&& \qquad{} +C\sup_{k}\ep(\varepsilon_{k}^{2}+Y^{2}_{0}+1)I\{
\widetilde
{Y}_{n-1}\leq
-\beta M\}
\end{eqnarray*}
and
\[
\ep Y^{2}_{n-1}I\{Y_{n-1}\leq
r,\widetilde{Y}_{n-1}\leq r\}+\ep
\widetilde{Y}^{2}_{n-1}I\{\widetilde{Y}_{n-1}\leq r\}\leq
C\sup_{k}\ep(\varepsilon_{k}^{2}+Y^{2}_{0}+1)I\{\widetilde
{Y}_{n-1}\leq
r\}.
\]
Since $\lim_{n\rightarrow\infty}\pr(\widetilde{Y}_{n}\leq x)=0$ for
any $x\in R$, we have
\[
\lim_{M\rightarrow\infty}\limsup_{n\rightarrow\infty}n^{-1}\sum
_{k=1}^{n}\widetilde{q}_{k}=0.
\]
This, together with $\ep\hat{\varepsilon}^{2}_{1}\rightarrow0$
as $K\rightarrow\infty$ and (\ref{xd1}), implies (\ref{p3}).
\end{pf}

\section{\texorpdfstring{Proof of Theorem \protect\ref{th3}}{Proof of Theorem 2.2}}\label{sec4}

To prove Theorem \ref{th3}, we need to establish the limiting
distribution for $\{Y_{t}\}$ as $t\rightarrow\infty$.

\begin{theorem}\label{th6} Let $\gamma=\delta=0$. Suppose either \textup{(H1)}
or \textup{(H2)} in Theorem \ref{th3} holds.
Then we have $\ep
|Y_{n}|=\mathrm{O}(\alpha^{n})$, $\sum_{t=0}^{\infty}I\{Y_{t}\leq r\}<\infty$
a.s. and $Y_{n}/\alpha^{n}\rightarrow\xi>0$ a.s., where $\xi$ is
defined in Theorem \ref{th3}.
\end{theorem}

From this theorem, we can see that
$\hat{\beta}_{n}-\beta\rightarrow Z$ a.s. for some random variable
$Z$. Thus,~$\hat{\beta}_{n}$ is not a strongly consistent
estimator for $\beta$. This explains why (\ref{eq1.5}) is the necessary
and sufficient condition for consistency of
$(\hat{\alpha}_{n},\hat{\beta}_{n})$. To prove Theorem \ref{th6}, we
need the following lemma.

\begin{lemma}\label{le3.2} Under the conditions of Theorem \ref{th6},
we have: \textup{(i)} $\ep|Y_{n}|=\mathrm{O}(\alpha^{n})$;
\textup{(ii)}~$\lim_{n\rightarrow\infty} Y_{n}/ n=\infty$ a.s. if
$\limsup_{n\rightarrow\infty}Y_{n}=\infty$ a.s.
\end{lemma}

\begin{pf} (i) Following the proof of Lemma \ref{le3.1},
we can prove that, under the conditions (H1) or~(H2) in Theorem
\ref{th3}, $\ep|Y_{n}|I\{|Y_{n}|\leq r\}=\mathrm{O}(1)$ if $\beta<1$ and
$\ep|Y_{n}|I\{|Y_{n}|\leq r\}=\mathrm{O}(n)$ if $\beta=1$. We next show that
$\ep|Y_{n}|=\mathrm{O}(\alpha^{n})$. Since
%
\begin{eqnarray}\label{p8}
\ep|Y_{n}|&=&\ep|Y_{n}|I\{Y_{n}\leq r\}\nonumber\\
&&{}+\sum_{k=0}^{n-1}\ep
|Y_{n}|I\{Y_{n}> r,Y_{n-1}> r,\ldots,Y_{k+1}> r,Y_{k}\leq r\}
\\
&&{}+\ep|Y_{n}|I\{Y_{n}> r,Y_{n-1}> r,\ldots,Y_{0}> r\}\nonumber
\end{eqnarray}
and $\ep|Y_{n}|I\{Y_{n}> r,Y_{n-1}> r,\ldots,Y_{0}>
r\}=\mathrm{O}(\alpha^{n})$, we only need to estimate the second term on the
right-hand side of (\ref{p8}). Set $B_{k}=\{Y_{n}> r,Y_{n-1}>
r,\ldots,Y_{k+1}> r,Y_{k}\leq r\}$. On $B_{k}$, we have
$Y_{n}=\sum_{j=k+2}^{n}\alpha^{n-j}\varepsilon_{j}+\alpha^{n-k-1}Y_{k+1}$.
Thus by noting that $\ep|Y_{n}|I\{B_{k}\}\leq\ep\max_{0\leq i\leq
n}
|\sum_{j=i+2}^{n}\alpha^{n-j}\varepsilon_{j}|I\{B_{k}\}+\alpha
^{n-k-1}\ep
|Y_{k+1}|I\{B_{k}\}$ and $\ep|Y_{n}|I\{Y_{n}\leq r\}=\mathrm{O}(n)$, we have
\[
\sum_{k=1}^{n-1}\ep|Y_{n}|I\{B_{k}\}\leq\ep\max_{0\leq i\leq n}
\Biggl|\sum_{j=i+2}^{n}\alpha^{n-j}\varepsilon_{j}\Biggr|+\mathrm{O}(1)\sum
_{k=0}^{n-1}\alpha
^{n-k-1}k=\mathrm{O}(\alpha^{n}).
\]
This together with (\ref{p8}) gives $\ep|Y_{n}|=\mathrm{O}(\alpha^{n})$.

(ii) For any $M>1$, define
$\mathbf{A}_{n}=\bigcup_{t=n}^{\infty}\{Y_{t}\leq t^{3/2}\}$. Let
$\delta>0$ and $T>\max(r,0)$ satisfy
$\alpha>1+\delta+8T^{-1/8}$. Define $\tau=\max\{k\dvtx Y_{-1}\leq
T,\ldots, Y_{k}\leq T, Y_{k+1}>T, k\geq-1\}$, $Y_{-1}=0$. We can see
that $\tau<\infty$ a.s. and $\{\tau=k\}=\{Y_{-1}\leq T,\ldots,
Y_{k}\leq T, Y_{k+1}>T\}$ is
$\sigma(Y_{0},\varepsilon_{1},\ldots,\varepsilon_{k+1})$ measurable.
For any $n_{0}+3<n$, $M>0$, $T>M$
\begin{eqnarray*}
\pr(\mathbf{A}_{n})&\leq&
\pr(\tau>n_{0})+\sum_{k=-1}^{n_{0}}\pr(\tau=k,
\mathbf{A}_{n})\\
&\leq& \pr(\tau>n_{0})+
\sum_{k=-1}^{n_{0}}\pr\Biggl(\tau=k,
\mathbf{A}_{n},\bigcap_{j=k+2}^{\infty}\bigl\{|\varepsilon
_{j}|\leq
\delta\bigl((j-k-2)^{2}+T\bigr)\bigr\}\Biggr)\\
& &{}+
\sum_{k=-1}^{n_{0}}\pr\Biggl(\tau=k,
\mathbf{A}_{n},\bigcup_{j=k+2}^{\infty}\bigl\{|\varepsilon_{j}|>
\delta\bigl((j-k-2)^{2}+T\bigr)\bigr\}\Biggr).
\end{eqnarray*}
Note that on the event
\[
\mathbf{B}:=\Biggl\{\tau=k,
\bigcap_{j=k+2}^{\infty}\bigl\{|\varepsilon_{j}|\leq
\delta\bigl((j-k-2)^{2}+T\bigr)\bigr\}\Biggr\},
\]
since $\alpha>1+\delta+8T^{-1/8}$ and
$\frac{(t-k-1)^{2}-(t-k-2)^{2}}{(t-k-2)^{2}+T}<8T^{-1/8}$ for $t\geq
k+3$, we have
\begin{eqnarray*}
Y_{k+1}&>&T>r,\\
Y_{k+2}&=&\alpha Y_{k+1}+\varepsilon_{k+2}\geq
\alpha T-\delta T>T+1,\\
 Y_{k+3}&=&\alpha
Y_{k+2}+\varepsilon_{k+3}\geq\alpha(T+1)-\delta(1+T)>T+2^{2},\\
\vdots&&\\
Y_{t}&=&\alpha
Y_{t-1}+\varepsilon_{t}>\alpha\bigl((t-k-2)^{2}+T\bigr)-\delta
\bigl((t-k-2)^{2}+T\bigr)\\
&>&(t-k-1)^{2}+T
\end{eqnarray*}
for any $t\geq k+1$. That is, for any $t$ satisfying
$t-n_{0}-1>t^{3/4}$ and $k\leq n_{0}$, we have $Y_{t}>t^{3/2}$ on
event $\mathbf{B}$. Thus, for $n$ satisfying $n-n_{0}-1>n^{3/4}$ and
$k\leq n_{0}$, we have
\[
\Biggl\{\tau=k,
\mathbf{A}_{n},\bigcap_{j=k+2}^{\infty}\bigl\{|\varepsilon
_{j}|\leq
\delta\bigl((j-k-2)^{2}+T\bigr)\bigr\}\Biggr\}=\varnothing
\]
and
\begin{eqnarray*}
\pr(\mathbf{A}_{n})&\leq& \pr(\tau>n_{0})+
\sum_{k=-1}^{n_{0}}\pr\Biggl(\tau=k,
\bigcup_{j=k+2}^{\infty}\bigl\{|\varepsilon_{j}|>
\delta\bigl((j-k-2)^{2}+T\bigr)\bigr\}\Biggr)\\
&=&\pr(\tau>n_{0})+
\sum_{k=-1}^{n_{0}}\pr(\tau=k)\pr\Biggl(
\bigcup_{j=k+2}^{\infty}\bigl\{|\varepsilon_{j}|>
\delta\bigl((j-k-2)^{2}+T\bigr)\bigr\}\Biggr)\\
&\leq&\pr(\tau>n_{0})+
\sum_{k=-1}^{n_{0}}\pr(\tau=k)\sum_{j=k+2}^{\infty
}\pr\bigl(
|\varepsilon_{j}|>\delta\bigl((j-k-2)^{2}+T\bigr)\bigr).
\end{eqnarray*}
Letting $n\rightarrow\infty$ and then $n_{0}\rightarrow\infty$
implies that for any $M>0$,
\begin{eqnarray*}
\pr\Biggl(\bigcap_{n=1}^{\infty}\mathbf{A}_{n}\Biggr) &\leq&
\delta^{-2}\sum_{k=-1}^{\infty}\pr(\tau=k)\sum
_{j=k+2}^{\infty}\frac{\ep|\varepsilon_{1}|^{2}}
{((j-k-2)^{2}+T)^{2}}\\
&\leq
&CT^{-1}\sum_{k=1}^{\infty}k^{-2}\rightarrow0\qquad  \mbox{as }
T\rightarrow\infty.
\end{eqnarray*}
The
lemma is proved.
\end{pf}

\begin{pf*}{Proof of Theorem \protect\ref{th6}}
$\ep|Y_{n}|=\mathrm{O}(\alpha^{n})$ follows from Lemma \ref{le3.2}(i). We next
give the proofs for the other conclusions.

\textit{Proof under} (H1). Define $X_{m}$ by the
equations $X_{n}=\varepsilon_{n}+X^{+}_{n-1}-\beta^{+} X_{n-1}^{-}$,
\mbox{$X_{0}=Y_{0}$}. Then we have $Y_{n}\geq X_{n}$ for any $n\geq0$. If
$\beta=1$, then $X_{n}=\sum_{k=1}^{n}\varepsilon_{k}+X_{0}$ and
$\limsup_{n\rightarrow\infty}Y_{n}=\limsup_{n\rightarrow\infty
}X_{n}=\infty$
a.s. If $\beta<1$, then by Theorem \ref{th4.1},
$X_{n}\rightarrow\infty$ in probability. So $Y_{n}\rightarrow\infty$
in probability, which implies
$\limsup_{n\rightarrow\infty}Y_{n}=\infty$ a.s. By~Lem\-ma~\ref{le3.2}
we have $Y_{n}/n\rightarrow\infty$ a.s. Hence
$\sum_{t=0}^{\infty}I\{Y_{t}\leq0\}<\infty$ a.s. Thus $\xi$ in
(\ref{newadd-11}) or (\ref{eq2.4}) is well defined and
$Y_{n}/\alpha^{n}\rightarrow\xi$ a.s.

Now we prove $\xi>0$ a.s. Let $e_{n}=\varepsilon_{n}-\beta
Y^{-}_{n-1}$. We have $\ep|e_{n}|=\mathrm{O}(n)$ for $\beta\leq1$. Define
$m=\sup\{n\dvtx \alpha^{n}/n<M\}$. Then $m\sim\log M/\log
\alpha\rightarrow\infty$ as $M\rightarrow\infty$. By the inequality
$(a+b)^{+}\geq a-|b|$, we have
%
%
\begin{eqnarray}\label{p9}
\frac{Y_{n}}{\alpha^{n}}=\frac{e_{n}}{\alpha^{n}}+\frac
{Y^{+}_{n-1}}{\alpha^{n-1}}
\geq
\frac{e_{n}}{\alpha^{n}}-\frac{|e_{n-1}|}{\alpha^{n-1}}+\frac
{Y^{+}_{n-2}}{\alpha^{n-2}}\geq
\cdots\geq
-\sum_{k=m+1}^{n}\frac{|e_{k}|}{\alpha^{k}}+\frac{Y^{+}_{m}}{\alpha^{m}}.
\end{eqnarray}
From (\ref{p9}) we can get $\xi\geq
-\sum_{k=m+1}^{\infty}\frac{|e_{k}|}{\alpha^{k}}+\frac
{Y^{+}_{m}}{\alpha^{m}}$
a.s. Since $M(m+1)/\alpha^{m+1}\leq1$, we have
\[
\sum_{k=m+1}^{\infty}(M\alpha^{-k})\ep|e_{k}|\leq
C\sum_{k=m+1}^{\infty}\frac{Mk}{\alpha^{k}}\leq
C\sum_{k=m+1}^{\infty}\frac{k-m}{\alpha^{k-m-1}}\frac{Mk}{\alpha
^{m+1}(k-m)}\leq
C\sum_{k=0}^{\infty}\frac{k+1}{\alpha^{k}},
\]
where $C$ does not depend on $M$.
So we have
%
%
\begin{eqnarray}\label{eq8}
\limsup_{M\rightarrow\infty}
\pr\Biggl(M\sum_{k=m+1}^{\infty}\frac{|e_{k}|}{\alpha^{k}}\geq
\eta\Biggr)\leq C\eta^{-1}\rightarrow0
\end{eqnarray}
as $\eta\rightarrow\infty$. It is easy to see that
$MY^{+}_{m}/\alpha^{m}\geq Y_{m}/m\rightarrow\infty$ a.s. as
$M\rightarrow\infty$. Hence, by (\ref{p9}) and~(\ref{eq8}),
$\pr(\xi\leq0)=\pr(M\xi\leq0)\leq\pr(Y_{m}/m\leq
\eta)+\pr(M\sum_{k=m+1}^{\infty}\frac
{|e_{k}|}{\alpha
^{k}}\geq
\eta)\rightarrow0 $ by letting $M\rightarrow\infty$ first and
then $\eta\rightarrow\infty$. This proves $\xi>0$ a.s.

\textit{Proof under} (H2). We first assume that
$\limsup_{n\rightarrow\infty} Y_{n}=\infty$ a.s. Then it follows
from Lemma~\ref{le3.2} that $Y_{n}/n\rightarrow\infty$ a.s. and
hence $\sum_{t=1}^{\infty}I\{Y_{t}\leq r\}<\infty$ a.s.,
$Y_{n}/\alpha^{n}\rightarrow\xi$ a.s. By writing
$Y_{n}=e_{n}+\alpha Y^{+}_{n-1}$, where $e_{n}=\varepsilon_{n}+\beta
Y_{n-1}I\{Y_{n-1}\leq r\}-\alpha Y_{n-1}I\{0\leq Y_{n-1}\leq r\}$ if
$r>0$, and $e_{n}=\varepsilon_{n}+\beta Y_{n-1}I\{Y_{n-1}\leq
r\}+\alpha Y_{n-1}I\{r<Y_{n-1}\leq0\}$ if $r\leq0$, we can show
that $\xi>0$ a.s. following the proof of Theorem \ref{th6} under
(H1).

It remains to show that $\limsup_{n\rightarrow\infty} Y_{n}=\infty$
a.s. We claim that if, for all $y\leq r$,
%
%
\begin{equation}\label{eq5}
\pr(Y_{t}< r \mbox{ for all $t\geq0$}|Y_{0}=y)=0,
\end{equation}
then $\limsup_{n\rightarrow\infty}Y_{n}=\infty$ a.s. The proof is
similar to that of \cite{13}. Let $c>|r|$. Since
for any $x>0$, $\pr(\varepsilon_{1}\leq x)<1$, we have for all
$r\leq y\leq c$,
\[
\pr(Y_{1}\geq c|Y_{0}=y)=\pr(\alpha
y+\varepsilon_{1}\geq c)\geq\pr\bigl(\varepsilon_{1}\geq c(1+\alpha)\bigr)>0,
\]
which yields that for any $c>|r|$
\[
\inf_{r\leq y\leq c}\pr(Y_{t}\geq c \mbox{ for some
$t>0$}|Y_{0}=y)>0.
\]
Then by Proposition 5.1 in \cite{11}, for any initial distribution
on $Y_{0}$,
%
%
\begin{equation}\label{eq7}
\{Y_{t}\in[r,c) \mbox{ infinitely often}\}\subseteq
\{Y_{t}\in[c,\infty) \mbox{ infinitely often}\}.
\end{equation}
Using similar arguments in \cite{13}, we can see
that if for all $y\in R$
%
%
\begin{equation}\label{eq4}
\pr(Y_{t}\geq r \mbox{ for some $t$}|Y_{0}=y)=1,
\end{equation}
then
\[
\pr(Y_{t}\geq r \mbox{ infinitely often})=1
\]
and hence by (\ref{eq7}) we have
$\pr(Y_{t}\in[c,\infty) \mbox{ infinitely often})=1$ for any
$c>0$. This yields $\limsup_{n\rightarrow\infty}Y_{n}=\infty$ a.s.

Now it suffices to show that (\ref{eq4}) or, equivalently, (\ref{eq5})
holds. Note that (\ref{eq5}) is a~direct consequence of the
following results:
%
%
\begin{equation}\label{eq6}
\lim_{n\rightarrow\infty}\pr\Biggl(\max_{1\leq k\leq
n}\Biggl(\sum_{i=1}^{k}\beta^{k-i}\varepsilon_{i}+\beta^{k}y
\Biggr)\leq
r\Biggr)=0 \qquad\mbox{for } y\leq r.
\end{equation}
If $\beta=1$, then (\ref{eq6}) holds by the law of iterated
logarithm. If $\beta\leq0$, we have
\begin{eqnarray*}
&&\Biggl\{\max_{1\leq k\leq
n}\Biggl(\sum_{i=1}^{k}\beta^{k-i}\varepsilon_{i}+\beta^{k}y
\Biggr)\leq
r\Biggr\}\\
&&\quad \subseteq
\Biggl\{\sum_{i=1}^{k}\beta^{k-i}\varepsilon_{i}+\beta
^{k}y=\varepsilon
_{k}+\beta
\Biggl(\sum_{i=1}^{k-1}\beta^{k-1-i}\varepsilon_{i}+\beta
^{k-1}y
\Biggr)\leq
r, 1\leq k\leq n\Biggr\}\\
&&\quad  \subseteq\{
\varepsilon_{k}\leq r+|\beta r|, 1\leq k\leq n\}.
\end{eqnarray*}
Therefore
\[
\pr\Biggl(\max_{1\leq k\leq
n}\Biggl(\sum_{i=1}^{k}\beta^{k-i}\varepsilon_{i}+\beta^{k}y
\Biggr)\leq
r\Biggr)\leq\pr\Bigl(\max_{1\leq k\leq n}\varepsilon_{k}\leq
r+|\beta r|\Bigr)\rightarrow0.
\]
It remains to prove (\ref{eq6}) for $0<\beta<1$. Set
$k_{j}=jn^{1/2}$ for $1\leq j\leq n^{1/2}$. Then for any $x>0$ we
have
\begin{eqnarray*}
\pr\Biggl(\max_{1\leq k\leq
n}\Biggl(\sum_{i=1}^{k}\beta^{k-i}\varepsilon_{i}\Biggr)\leq
x\Biggr)&\leq& \pr\Biggl(\max_{1\leq j\leq
n^{1/2}}\Biggl(\sum_{i=1}^{k_{j}}\beta^{k_{j}-i}\varepsilon_{i}
\Biggr)\leq
x\Biggr)\\
&\leq& \pr\Biggl(\max_{1\leq j\leq
n^{1/2}}\Biggl(\sum_{i=k_{j}-n^{1/4}}^{k_{j}}\beta
^{k_{j}-i}\varepsilon
_{i}\Biggr)\leq
2x\Biggr)\\
& &{}+
\pr\Biggl(\max_{1\leq j\leq
n^{1/2}}\Biggl(\sum_{i=1}^{k_{j}-n^{1/4}-1}\beta
^{k_{j}-i}\varepsilon
_{i}\Biggr)\geq
x\Biggr).
\end{eqnarray*}
Since $\ep|\varepsilon_{1}|<\infty$, we have
\[
\pr\Biggl(\max_{1\leq j\leq
n^{1/2}}\Biggl(\sum_{i=1}^{k_{j}-n^{1/4}}\beta^{k_{j}-i}\varepsilon
_{i}\Biggr)\geq
x\Biggr)\leq Cn^{1/2}\sum_{j=n^{1/4}}^{\infty}\beta^{j}\rightarrow
0.
\]
By independence, we have
\[
\pr\Biggl(\max_{1\leq j\leq
n^{1/2}}\Biggl(\sum_{i=k_{j}-n^{1/4}}^{k_{j}}\beta
^{k_{j}-i}\varepsilon
_{i}\Biggr)\leq
2x\Biggr)=\Biggl(\pr\Biggl(\sum_{j=1}^{n^{1/4}+1}\beta
^{j-1}\varepsilon
_{j}\leq
2x\Biggr)\Biggr)^{n^{1/2}}.
\]
Also
\[
\pr\Biggl(\sum_{j=1}^{n^{1/4}+1}\beta^{j-1}\varepsilon_{j}\leq
2x\Biggr) \leq
\pr\Biggl(\sum_{j=1}^{\infty}\beta^{j-1}\varepsilon_{j}\leq
3x\Biggr)+\Delta_{n},
\]
where $\Delta_{n}\leq C\sum_{j=n^{1/4}}^{\infty}\beta
^{j}\rightarrow
0$ as $n\rightarrow\infty$. So it suffices to show that for any
$x>0$, $\pr(\sum_{j=1}^{\infty}\beta^{j-1}\varepsilon
_{j}\leq
x)<1$. In fact, if there exists some $x>0$ such that
\[
1=\pr\Biggl(\sum_{j=1}^{\infty}\beta^{j-1}\varepsilon_{j}\leq
x\Biggr)=\ep F\Biggl(x-\sum_{j=2}^{\infty}\beta^{j-1}\varepsilon_{j}\Biggr),
\]
where $F(\cdot)$ is the distribution function of $\varepsilon_{1}$,
then $F(x-\sum_{j=2}^{\infty}\beta^{j-1}\varepsilon_{j})=1$ a.s.
That is, $\sum_{j=2}^{\infty}\beta^{j-1}\varepsilon_{j}=-\infty$ a.s.
This is impossible since $0<\beta<1$ and
$\ep|\varepsilon_{1}|<\infty$.
\end{pf*}

We are now ready to prove Theorem \ref{th3}.

\begin{pf*}{Proof of Theorem \protect\ref{th3}} Note that
\[
\frac{\alpha^{n}(\hat{\alpha}_{n}-\alpha)}{\alpha^{2}-1}=\frac
{\alpha
^{-n}\sum_{t=1}^{n-1}I(Y_{t}>
r)Y_{t}\varepsilon_{t+1}}{\alpha^{-2n}(\alpha^{2}-1)\sum_{t=1}^{n-1}I(Y_{t}>
r)Y^{2}_{t}}.
\]
Since
$Y_{n}/\alpha^{n}\rightarrow\xi>0$ a.s., we have
$(\alpha^{2}-1)\alpha^{-2n}\sum_{t=1}^{n-1}Y_{t}^2I\{Y_{t}>r\}
\rightarrow
\xi^{2}$ a.s. By the fact $\ep|Y_{n}|I\{Y_{n}\leq r\}=\mathrm{O}(n)$ for
$\beta\leq1$, we have
\[
\alpha^{-n}\Biggl(\sum_{t=1}^{n-1}Y_{t}\varepsilon_{t+1}-\sum
_{t=1}^{n-1}Y_{t}\varepsilon_{t+1}I\{Y_{t}>r\}\Biggr)
\rightarrow0\qquad  \mbox{a.s.}
\]
We next prove that
$\alpha^{-n}(\sum_{t=1}^{n-1}Y_{t}\varepsilon_{t+1}-\xi\sum
_{t=1}^{n-1}\alpha^{t}\varepsilon_{t+1})\rightarrow
0$ in probability. For $K>0$, let
\[
\widetilde{\varepsilon}_{t}=\varepsilon_{t}I\{|\varepsilon_{t}|\leq
K\},\qquad 1\leq t\leq n.
\]
We have
$\alpha^{-n}\ep{|}\sum_{t=1}^{n-1}Y_{t}(\varepsilon
_{t+1}-\widetilde
{\varepsilon}_{t+1}){|}\leq
C\ep|\varepsilon_{0}|I\{|\varepsilon_{0}|> K\}\rightarrow0 $ as
$K\rightarrow\infty$. So it suffices to prove that
$\alpha^{-n}(\sum_{t=1}^{n-1}Y_{t}\widetilde{\varepsilon}_{t+1}
-\xi\sum_{t=1}^{n-1}\alpha^{t}\widetilde{\varepsilon}_{t+1}
)\rightarrow
0$ in probability, which follows from
$Y_{n}/\alpha^{n}\rightarrow\xi$ a.s. and
$|\widetilde{\varepsilon}_{t}|\leq K$. Hence
\[
\frac{\alpha^{n}(\hat{\alpha}_{n}-\alpha)}{\alpha^{2}-1}-\frac
{\alpha
^{-n}\sum_{t=1}^{n-1}\alpha^{t}\varepsilon_{t+1}}{\xi}\rightarrow
0\qquad \mbox{in probability}.
\]
Note that $Y_{[n/2]}/\alpha^{[n/2]}\rightarrow\xi$ a.s. and
$\alpha^{-n}(\sum_{t=1}^{n-1}\alpha^{t}\varepsilon
_{t+1}-\sum
_{t=[n/2]+1}^{n-1}\alpha^{t}\varepsilon_{t+1})\rightarrow
0$ in probability. We have
\[
\frac{\alpha^{n}(\hat{\alpha}_{n}-\alpha)}{\alpha^{2}-1}-\frac
{\alpha
^{-n}\sum_{t=[n/2]+1}^{n-1}\alpha^{t}\varepsilon_{t+1}}
{Y_{[n/2]}/\alpha^{[n/2]}}\rightarrow0 \qquad\mbox{in probability}.
\]
By the independence between
$\sum_{t=[n/2]+1}^{n-1}\alpha^{t}\varepsilon_{t+1}$ and $Y_{[n/2]}$,
we see that
\[
\Biggl(\alpha^{-n}\sum_{t=[n/2]+1}^{n-1}\alpha^{t}\varepsilon_{t+1},
Y_{[n/2]}/\alpha^{[n/2]}\Biggr)\Rightarrow(\eta^{*},\xi^{*}),
\]
which finishes the proof.
\end{pf*}

\section{A further result when $\alpha\beta=1$}\label{sec5}

We next consider the LSE of $(\alpha,\beta)$ under the constraints
$\alpha\beta=1$. We estimate $\alpha$ by minimizing $Q_{n}(x)$,
where
\[
Q_{n}(x)=\sum_{t=2}^{n}(Y_{t}-xY_{t-1}I\{Y_{t-1}<
r\}-x^{-1}Y_{t-1}I\{Y_{t-1}\geq r\})^{2}.
\]
Pham, Chan and Tong \cite{13} showed that the estimator
$\hat{\alpha}_{n}$, by minimizing $Q_{n}(x)$ under $\alpha\beta=1$
and $\alpha<0$, is strongly consistent. The following theorem shows
that $\hat{\alpha}_{n}$ is still strongly consistent under
$\alpha\beta=1$ and $\alpha>0$.

\begin{theorem}\label{th4} Let $\gamma=\delta=0$, $\alpha\beta=1$ and
$0<\alpha\neq1$. Assume that $\pr(\varepsilon_{1}\leq x)<1$ and
$\pr(\varepsilon_{1}\geq x)<1$ for any $x\in R$. Then
$\hat{\alpha}_{n}$ obtained by minimizing $Q_{n}(x)$ is strongly
consistent.
\end{theorem}

\begin{pf} We only prove the
theorem for $\alpha>1$. The proof for the other case $0<\alpha<1$ is
similar. We have
\begin{eqnarray*}
&&Q_{n}(x)-Q_{n}(\alpha)\\
&&\quad=(x-\alpha)^{2}\sum_{t=2}^{n}Y^{2}_{t-1}I\{Y_{t-1}>
r\}-2(x-\alpha)\sum_{t=2}^{n}\varepsilon_{t}Y_{t-1}I\{Y_{t-1}>
r\}\\
&&\qquad{}+(x^{-1}-\alpha^{-1})^{2}\sum_{t=2}^{n}Y^{2}_{t-1}I\{Y_{t-1}\leq
r\}\\
&&\qquad{}-2(x^{-1}-\alpha^{-1})\sum_{t=2}^{n}\varepsilon_{t}Y_{t-1}I\{
Y_{t-1}\leq
r\}\\
&&\quad\geq (x-\alpha)^{2}\sum_{t=2}^{n}Y^{2}_{t-1}I\{Y_{t-1}>
r\}-2(x-\alpha)\sum_{t=2}^{n}\varepsilon_{t}Y_{t-1}I\{Y_{t-1}>
r\}-\sum_{t=2}^{n}\varepsilon^{2}_{t}.
\end{eqnarray*}
By Theorem \ref{th6}, we can see that
\begin{eqnarray*}
\frac{1}{\alpha^{2n}}\sum_{t=2}^{n}Y^{2}_{t-1}I\{Y_{t-1}>r\}
&\rightarrow&
(\alpha^{2}-1)^{-1}\xi  \qquad\mbox{a.s.};\\
\sum_{t=2}^{n}\varepsilon_{t}Y_{t-1}I\{Y_{t-1}>
r\}&=&\mathrm{O}(\alpha^{3n/2})  \qquad\mbox{a.s};\\
\sum_{t=2}^{n}Y^{2}_{t-1}I\{Y_{t-1}\leq r\}&=&\mathrm{O}(1)
\qquad\mbox{a.s.};\\
\sum_{t=2}^{n}\varepsilon_{t}Y_{t-1}I\{Y_{t-1}\leq r\}&=&\mathrm{O}(1)  \qquad\mbox{a.s.}
\end{eqnarray*}
Hence for any $\delta>0$, we have
\[
\lim_{n\rightarrow\infty}\inf_{x:|x-\alpha|>\delta}
\bigl(Q_{n}(x)-Q_{n}(\alpha)\bigr)=\infty \qquad \mbox{a.s.}
\]
Since $Q_{n}(x)$ is continuous on $[\alpha-\delta,\alpha+\delta]$,
it always admits a minimum on this interval. This shows that
$\limsup_{n\rightarrow\infty}|\hat{\alpha}_{n}-\alpha|\leq\delta$
a.s. for any $\delta>0$ and completes the proof.
\end{pf}

\section*{Acknowledgements}
The authors thank the referee and the editor
Holger Rootz\'{e}n for their helpful comments.
Shiqing Ling's research supported
by HK RGC CERG 601607 and CERG 602609.
Qi-Man Shao's research partially
supported by Hong Kong RGC CERG 602608.

\printhistory

\end{document}